\magnification=\magstep1
\vsize=23.1true cm
\hsize=16.2 true cm
\nopagenumbers
\topskip=1true cm
\headline={\tenrm\hfil\folio\hfil}
\raggedbottom
\abovedisplayskip=3mm
\belowdisplayskip=3mm
\abovedisplayshortskip=0mm
\belowdisplayshortskip=2mm
\normalbaselineskip=12pt
\normalbaselines
\font\b=cmssdc10 scaled\magstep1

\centerline{\bf A CERTAIN $p$-ADIC SPECTRAL THEOREM}
\bigskip
\centerline{01/13/07}
\bigskip
\bigskip
\centerline{\it R. L. Baker}
\centerline{\it University of Iowa}
\centerline{\it Iowa City, Iowa 52242}
\vfill\eject
\centerline{\bf A CERTAIN $p$-ADIC SPECTRAL THEOREM}
\bigskip
\centerline{\it R. L. Baker}
\centerline{\it  University of Iowa}
\centerline{\it  Iowa City, Iowa 52242}
\vfill\eject
\noindent
{\bf ABSTRACT.} We extend a $p$-adic spectral theorem of M. M. Vishik to a certain class of $p$-adic Banach algebras.
This class includes inductive limits of finite-dimensional $p$-Banach algebras of the form $B({\cal X})$, where $\cal X$
is a $p$-adic Banach space of the form ${\cal X}\simeq\Omega_p(J)$, $J$ being a finite nonempty set. In particular,
we present a $p$-adic spectral theorem for $p$-adic UHF algebras and $p$-adic TUHF algebras (Triangular UHF Algebras).

\vfill\eject
\parindent=1cm
\centerline{\bf 1. INTRODUCTION}
\vskip 1cm
In the article, {\it On certain Banach Limits of Triangular Matrix Algebras} ({\bf [B1]}) the author of the present paper
proved, using purely Banach-algebraic techniques, that the supernatural number associated to an arbitrary
{\bf triangular UHF (TUHF)} Banach algebra over the complex number field {\b C} is an invariant of the algebra, provided
that the algebra satisfies certain ``local dimensionality conditions.'' The proof, although relying only on classical 
complex  Banach algebra techniques, makes essential use of the classical spectral theorem (the Riesz functional calculus)
for complex Banach algebras. For any prime number $p$, let $\Omega_p$ be the $p$-adic counterpart of {\b C}.
Replacing {\b C} by $\Omega_p$ in the definition of complex TUHF Banach algebras, we obtain the definition of $p$-adic
$p$-adic TUHF algebras. The results in {\bf [B1]} for complex TUHF algebras can be 
duplicated for $p$-adic TUHF algebras, provided that a sufficiently general $p$-adic version of the Riesz functional
calculus can be developed for $p$-adic Banach algebras ({\bf [B2]}).
The purpose to the present paper is to prove a $p$-version of the classical spectral theorem for complex Banach algebras
and to apply this version to $p$-adic TUHF algebras. The $p$-adic spectral theorem that we develop is an extension of 
a certain $p$-adic spectral theorem of M. M. Vishik ({\bf [K]}, p. 149).

In the present section of the paper we put forth the preliminary material on $p$-adic analysis and $p$-adic Banach 
algebras that is necessary for proving an extended version of Vishik's Spectral Theorem (Theorem 1.18). Then in Section 2
we prove the extended version of Vishik's Spectral Theorem (Theorem 2.10). Finally, in Section 3, we show that 
this extended spectral theorem provides a {\bf Spectral Theorem for $p$-adic TUHF and UHF Algebras}, which is Theorem 3.3.

Let $p$ be a prime number, and let {\b Q} be the field of rational numbers. Let $|\cdot|_p$ be the function defined on 
{\b Q} by
$$\bigg|{a\over b}\bigg|_p=p^{{\rm ord}_p b-{\rm ord}_p a},\quad |0|_p=0.$$
Here ${\rm ord}_p$ of a non-negative integer is the highest power of $p$ dividing the integer. Then $|\cdot|_p$ 
is a norm on {\b Q}. The field {\b Q}$_p$ is
defined to be the completion of {\b Q} under the norm $|\cdot|_p$. Unlike the case of the real numbers {\b R}, 
whose algebraic closure {\b C} is only a
quadratic extension of {\b R}, the algebraic closure $\overline{\hbox{\b Q}}_p$ of {\b Q}$_p$ 
has infinite degree over {\b Q}$_p$. However, the norm
 $|\cdot|_p$ on {\b Q}$_p$ can be extended to a norm $|\cdot|_p$ on $\overline{\hbox{\b Q}}_p$. But it turns out 
 that  $\overline{\hbox{\b Q}}_p$ is not
complete under this extended norm. Thus, in order to do analysis, we must take a larger field than 
$\overline{\hbox{\b Q}}_p$. We denote the completion of $\overline{\hbox{\b Q}}_p$ under the norm  
$|\cdot|_p$ by $\Omega_p$, that is,
$$\Omega_p=\buildrel{\wedge\;\;}\over{\overline{\hbox{\b Q}}_p},$$
where $^\wedge$ means completion with respect to  $|\cdot|_p$. ({\it Note}: The symbol ``{\b C}$_p$'' 
is sometimes used to denote $\Omega_p$ (see {\bf [K]}: p. 13).
\medskip
\proclaim 1.1 Definition. Let $r\ge 0$ be an non-negative real number, and let $a\in\Omega_p$ 
and let $\sigma\subseteq\Omega_p$. We have the following definitions.
$$\eqalign{D_a(r)&=\{\,x\in\Omega_p\;\;|\quad |x-a|_p\le r\,\};\cr
           D_a(r^-)&=\{\,x\in\Omega_p\;\;|\quad |x-a|_p<r\,\};\cr
	              D_\sigma(r)&=\{\,x\in\Omega_p\;\;|\quad{\rm dist}(x,\sigma)\le r\,\};\cr
		                 D_\sigma(r^-)&=\{\,x\in\Omega_p\;\;|\quad{\rm dist}(x,\sigma)<r\,\}.\cr}$$
If $b\in D_a(r)$, then $D_b(r)=D_a(r)$, and if $b\in D_a(r^-)$, then $D_b(r^-)=D_a(r^-)$. Thus, any point in a 
disc is its center. Hence $D_a(r), D_a(r^-)$
are both open and closed in the topological sense.
\medskip

\proclaim 1.2 Lemma. Let $\emptyset\not=\sigma\subseteq\Omega_p$ be compact. Then for any $s>0$ there exists 
$0<r\in |\Omega_p|_p$ and
$a_1,...,a_N\in\sigma$ such that $r<s$ and 
$$\eqalign{ D_\sigma(r)&=\bigcup\limits^N_{i=1}D_{a_i}(r);\cr
             \sigma&\subseteq\bigcup\limits^N_{i=1}D_{a_i}(r^-),\cr}$$
where the $D_{a_i}(r)$ are disjoint.
\medskip

\noindent
{\bf Proof.} Since $\sigma$ is compact, there exist $a_1,...,a_N\in\sigma$ such that
$$\sigma\subseteq\bigcup\limits^N_{i=1}D_{a_i}(s^-),\eqno (1)$$
where the union is disjoint. Now, for each $1\le i\le N$, $D_{a_i}(s^-)$ is closed, and hence
$\sigma\cap D_{a_i}(s^-)$ is compact. Therefore there exists $b_1,...,b_N\in\sigma$ such that for $1\le i\le N$,
$$|b_i-a_i|_p=\sup\{\,|x-a_i|_p\,:\,x\in\sigma\cap D_{a_i}(s^-)\,\}.$$
Because each $b_i\in D_{a_i}(s^-)$, we have, for $1\le i\le N$, $s>|b_i-a_i|_p$.
Now select any $r\in |\Omega_p|_p$ such that $r>0$ and
$$s>r>|b_1-a_1|_p,\ldots,|b_N-a_N|_p.$$
We claim that
$$\sigma\subseteq\bigcup\limits^N_{i=1}D_{a_i}(r^-),\eqno (2)$$
where $D_{a_1}(r),\ldots,D_{a_N}(r)$ are disjoint. To prove the claim, let $a\in\sigma$; then by (2),
$a\in D_{a_i}(s^-)$ for some $1\le i\le N$, and hence $a\in\sigma\cap D_{a_i}(s^-)$, which implies that
$$|a-a_i|_p\le \sup\{\,|x-a_i|_p\,:\,x\in\sigma\cap D_{a_i}(s^-)\,\}=|b_i-a_i|_p<r.$$
This shows that $a\in D_{a_i}(r^-)$, which proves (2). Now let $1\le i,j\le N$, with $i\not= j$. Suppose that
$D_{a_i}(r)\cap D_{a_j}(r)\not=\emptyset$, and let $a\in D_{a_i}(r)\cap D_{a_j}(r)$. Then we have
$$|a_i-a_j|_p\le\max\{\,|a_i-a|_p,|a-a_j|_p\,\}\le r<s,$$
hence $a_i\in D_{a_i}(s^-)\cap D_{a_j}(s^-)=\emptyset$. This contradiction proves the claim.
We now claim that
$$D_\sigma(r)=\bigcup\limits^N_{i=1}D_{a_i}(r). \eqno (3)$$
It is clear that $\bigcup\limits^N_{i=1}D_{a_i}(r)\subseteq D_\sigma(r)$. Let $x\in D_\sigma(r)$, then because $\sigma$
is compact, there exists $a\in\sigma$ such that $r\ge {\rm dist}(x,\sigma)=|x-a|_p$. Hence, $x\in D_a(r)$. Now,
$a\in\sigma$, so by (2), there exists $1\le j\le N$ such that $a\in D_{a_j}(r^-)\subseteq D_{a_j}(r)$. Therefore,
$D_a(r)=D_{a_j}(r)$, which shows that $x\in\bigcup\limits^N_{i=1}D_{a_i}(r)$. This proves (3).
\vrule height 6pt width 5pt depth 4pt

\medskip
\proclaim 1.3 Lemma. Let $\emptyset\not=\sigma\subseteq\Omega_p$. 
Let $r>0$ and let $b_1,\ldots,b_M$ be in $\Omega_p$, with
$$\sigma\subseteq\bigcup\limits^M_{i=1}D_{b_i}(r^-),\eqno (1)$$
where the $D_{b_i}(r)$ are disjoint. Then there exist $a_1,\ldots,a_N$ in $\sigma$ and
$\emptyset\not=I\subseteq\{\,1,\ldots,M\,\}$ such that the $D_{a_i}(r)$ are disjoint and
$$\eqalignno{D_\sigma(r^-)&=\bigcup\limits^N_{i=1}D_{a_i}(r^-)=\bigcup\limits_{i\in I}D_{b_i}(r^-);& (2)\cr
             D_\sigma(r)&=\bigcup\limits^N_{i=1}D_{a_i}(r)=\bigcup\limits_{i\in I}D_{b_i}(r).}$$
\medskip

\noindent
{\bf Proof.} Let $\emptyset\not=I\subseteq\{\,1,\ldots,M\,\}$ be the subset of $\{\,1,\ldots,M\,\}$ consisting of those
$1\le i\le M$ for which $\sigma\cap D_{b_i}(r^-)\not=\emptyset.$ Write
$$I=\{\,k_1,\ldots,k_N\,\}.$$
For each $1\le i\le N$, select an $a_i\in \sigma\cap D_{b_{k_i}}(r^-)$. Then for $1\le i\le N$, we have
$D_{a_i}(r^-)= D_{b_{k_i}}(r^-)$, and hence, by (1),
$$\sigma\subseteq \bigcup\limits^N_{i=1}D_{a_i}(r^-)=\bigcup\limits_{i\in I}D_{b_i}(r^-).\eqno (3)$$
The assumption that the $D_{b_i}(r)$ are disjoint easily implies that the $D_{a_i}(r)$ are disjoint. We claim that
$$\eqalignno{D_\sigma(r^-)&\subseteq \bigcup\limits^N_{i=1}D_{a_i}(r^-); & (4)\cr
D_\sigma(r)&\subseteq \bigcup\limits^N_{i=1}D_{a_i}(r).}$$
To prove (4) let $x\in D_\sigma(r^-)$. Then, by definition, ${\rm dist}(x,\sigma)<r$,
hence there exists $a\in\sigma$ such that $|x-a|_p<r$. Then (3) implies that there exist $1\le j\le N$ such that
$|a-a_j|_p<r$. Therefore, we have
$$|x-a_j|_p=|(x-a)+(a-a_j)|_p\le \max\{\,|x-a|_p,|a-a_j|_p\,\}<r.$$
Thus, $x\in D_{a_j}(r^-)$. This proves the first part (4). To prove the second part of (4), suppose that 
there exists $x\in D_\sigma(r)$ with
$x\not\in \bigcup\limits^N_{i=1}D_{a_i}(r)$. Then $|x-a_i|_p>r$ for $1\le i\le N$, hence we can find $0<r<s$ such that
$|x-a_i|_p>s$ for $1\le i\le N$. Now let $a\in\sigma$ be arbitrary, then by (3), $|a-a_i|_p<r$ for some $1\le i\le N$.
We then have
$$|x-a|_p=|(x-a_i)+(a_i-a)|_p=|x-a_i|_p>s.\eqno (5)$$
In (5) we have used the easily proved fact that for $y,z\in\Omega_p$, if $|y|_p>|z|_p$, then
$|y+z|_p=|y|_p$ (see {\bf [BBN]}: Theorem 2, p. 5). But $a\in\sigma$ is arbitrary, hence (5) implies
that ${\rm dist}(x,\sigma)\ge s>r$, which contradicts $x\in D_\sigma(r^-)$. This proves the second part of (4). 
It is easy to see that 
$$\eqalignno{\bigcup\limits^N_{i=1}D_{a_i}(r^-)&\subseteq D_\sigma(r^-); & (6)\cr
             \bigcup\limits^N_{i=1}D_{a_i}(r)&\subseteq  D_\sigma(r).}$$
Hence by (4) and (6), we see that (2) holds. \vrule height 6pt width 5pt depth 4pt

\medskip
\proclaim 1.4 Definition. Let $a\in\Omega_p$, and let $r>0$. 
A function $f:D_a(r)\to\Omega_p$ is said to be {\bf Krasner analytic} on
$D_a(r)$ iff $f$ can be represented by a power series on $D_a(r)$ of the form $f(x)=\sum^\infty_{k=0}c_k(x-a)^k$, where
$\lim\limits_{k\to\infty} r^k|c_k|_p=0$

\medskip
\proclaim 1.5 Definition. Let $\emptyset\not=\sigma\subseteq\Omega_p$, 
and let $r>0$. Define $B_r(\sigma)$ to be the set of all functions
$$f:D_\sigma(r)\to\Omega_p$$
such that $f$ is Krasner analytic on $D_a(r)$, whenever $a\in\Omega_p$ and $D_a(r)\subseteq D_\sigma(r)$. If $f$ is 
bounded on $D_\sigma(r)$, we write $\|f\|_r=\max\limits_{x\in D_\sigma(r)}|f(x)|_p$.

\medskip 
\proclaim 1.6 Definition. Let $\emptyset\not=\sigma\subseteq\Omega_p$. Define $L(\sigma)$ to be the set of all
$\Omega_p$-valued functions $f$ for which there exist $a_1,\ldots,a_N$ in $\Omega_p$ and $0<r\in |\Omega_p|_p$ such that
$$\sigma\subseteq\bigcup\limits^N_{i=1}D_{a_i}(r^-),$$
where the $D_{a_i}(r)$ are disjoint and $f$ is Krasner analytic on each $D_{a_i}(r)$. $L(\sigma)$ is called the set
of {\bf locally analytic functions} on $\sigma$

\medskip
\proclaim 1.7 Lemma. Let  $\emptyset\not=\sigma\subseteq\Omega_p$ be compact. Then
$$L(\sigma)=\bigcup\limits_{r>0}B_r(\sigma). \eqno (1)$$
\medskip

\noindent
{\bf Proof.} Let $f\in L(\sigma)$. Then by Definition 1.6 there exist $a_1,\ldots,a_N$ in
$\Omega_p$ and $0<s\in |\Omega_p|_p$ such that
$$\sigma\subseteq\bigcup\limits^N_{i=1}D_{a_i}(s^-), \eqno (2)$$
where the $D_{a_i}(s)$ are disjoint and $f$ is Krasner analytic on each $D_{a_i}(s)$. Then by (2) and Lemma 1.3,
there exists $\emptyset\not=I\subseteq \{\,1,\ldots,N\,\}$ such that
$$D_\sigma(r)=\bigcup\limits_{i\in I}D_{a_i}(s). \eqno (3)$$
Let $b\in\Omega_p$, with $D_b(s)\subseteq D_\sigma(s)$. Then by (3), there exists $i\in I$ such that
$b\in D_{a_i}(s)$, and consequently, $D_b(s)=D_{a_i}(s)$. Because $f$ is Krasner analytic on $D_{a_i}(s)$, we see that
$f$ is Krasner analytic on $D_b(s)$. This shows that $f\in B_s(\sigma)$. Because $f\in L(\sigma)$ is arbitrary, we get
$$L(\sigma)\subseteq\bigcup\limits_{r>0}B_r(\sigma). \eqno (4)$$
Now let $f\in \bigcup\limits_{r>0}B_r(\sigma)$, say, $f\in B_s(\sigma)$, $s>0$. By Lemma 1.3 there exists
$b_1,\ldots,b_M$ in $\Omega_p$ such that
$$\sigma\subseteq D_\sigma(s)=\bigcup\limits^M_{i=1}D_{b_i}(s). \eqno (5)$$
By Lemma 1.2 there exist $a_1,\ldots,a_N$ in $\sigma$ and $0<r\in |\Omega_p|_p$ such that $r<s$, the $D_{a_i}(r)$ are
disjoint, and
$$D_\sigma(r)=\bigcup\limits^N_{i=1}D_{a_i}(r),\quad \sigma\subseteq\bigcup\limits^N_{i=1}D_{a_i}(r^-). \eqno (6)$$
Let $1\le i\le N$ be arbitrary. Because $a_i\in\sigma$, (5) implies that there exists $1\le j\le M$ such that
$a_i\in D_{b_j}(s)$, which in tern implies that $D_{a_i}(s)=D_{b_j}(s)$. Now, $f\in B_s(\sigma)$, therefore $f$ is
Krasner analytic on $D_{b_j}(s)=D_{a_i}(s)$. Hence, by Definition 1.6., $f$ can be represented by a power series
on $D_{a_i}(s)$ of the form $f(x)=\sum\limits^\infty_{k=0}c_k(x-a_i)^k$, where
$\lim\limits_{k\to\infty}s^k|c_k|_p=0$. Since $r<s$, we have $D_{a_i}(r)\subseteq D_{a_i}(s)$, and hence, for
$x\in D_{a_i}(r)$, we have $f(x)=\sum\limits^\infty_{k=0}c_k(x-a_i)^k$, with $\lim\limits_{k\to\infty}r^k|c_k|_p=0$.
This shows that $f$ is Krasner analytic on $D_{a_i}(r)$. This proves that
$$\bigcup\limits_{r>0}B_r(\sigma)\subseteq L(\sigma). \eqno (7)$$
Statements (4) and (7) together imply (1). \vrule height 6pt width 5pt depth 4pt

In {\bf [K]}, p. 136, Koblitz defines the {\it locally analytic functions} on $\emptyset\not=\sigma\subseteq\Omega_p$
compact to be $L(\sigma)=\bigcup\limits_{r>0}B_r(\sigma)$. Thus, Lemma 1.7 shows that Definition 1.6 generalizes
Koblitz's locally analytic functions on $\sigma$ to the case where $\emptyset\not=\sigma\subseteq\Omega_p$
is not necessarily compact.

\medskip
\proclaim 1.8 Definition (The Shnirelman Integral). Let $r>0$, with $r\in |\Omega_p|_p$. Let $a\in\Omega_p$, and let
$f$ be an $\Omega_p-$valued function whose domain contains all $x\in\Omega_p$ such that $|\,x-a\,|_p=r$. Let
$\Gamma\in \Omega_p$, with $|\Gamma|_p=r$. Then the {\bf Shnirelman integral} of $f$ over the circle
$$\{\,x\in\Omega_p\,:\,|x-a|_p=r\,\}$$
is defined to be the following limit, provided the limit exists.
$$\int\limits_{a,\Gamma}f(x)\,dx=\lim\limits_{\buildrel{n\to\infty}\over{p\not\;\mid n}}{1\over n}
\sum\limits_{\xi^n=1}f(a+\xi\Gamma).$$

\medskip
\proclaim 1.9 Definition. For $\sigma\subseteq\Omega_p$ be compact, let $H_0(\overline{\sigma})$
denote the set of functions
$\varphi:\overline{\sigma}\to\Omega_p$ which are {\bf Krasner analytic} on $\overline{\sigma}$, i.e.,
\item\item{(1)} $\varphi$ is the limit of rational functions whose poles are contained in $\sigma$, the limit being
uniform in any set of the form
$$\overline{D}_\sigma(r)=\{\,z\in\Omega_p\;|\;{\rm dist}(z,\sigma)\ge r\,\},\quad \sigma\subseteq D_\sigma(r^-), r>0.$$
\item\item{(2)} $\lim\limits_{|z|_p\to\infty}\varphi(z)=0$.

\medskip
\proclaim 1.10 Lemma. Let $\emptyset\not=\sigma\subseteq\Omega_p$ be compact. Let $0<r_1\le r_2$,
and let $\Gamma_1,\Gamma_2$ be in $\Omega_p$, with $|\Gamma_1|_p=r_1$, $|\Gamma_2|_p=r_2$. Assume that
$$a_1,\ldots,a_M;b_1,\ldots,b_N\in\Omega_p$$
are given, with
$$\eqalign{D_\sigma(r_1)&=\bigcup\limits^M_{i=1}D_{a_i}(r_1),\quad \sigma \subseteq
                          \bigcup\limits^M_{i=1}D_{a_i}({r_1}^-);\cr
           D_\sigma(r_2)&=\bigcup\limits^N_{i=1}D_{b_i}(r_2), \quad \sigma \subseteq
                          \bigcup\limits^N_{i=1}D_{b_i}({r_2}^-),\cr}$$
where the $D_{a_i}(r_1)$ and the $D_{b_i}(r_2)$ are disjoint.
Then for $\phi\in H_0(\overline{\sigma})$ and $f\in B_{r_2}(\sigma)$, we have $f\in B_{r_1}(\sigma)$, and the following
sums exist and are equal.
$$\sum\limits^M_{i=1}\,\int\limits_{a_i,\Gamma_1}\,f(x)(x-a_i)\varphi(x)\,dx=
  \sum\limits^N_{i=1}\,\int\limits_{b_i,\Gamma_2}\,f(x)(x-b_i)\varphi(x)\,dx.$$

\noindent
{\bf Proof.} See {\bf [K]}:Lemma 8, p. 138.

\medskip
\proclaim 1.11 Lemma. Let $a\in\Omega_p$, and let $0<r$. Let $f:D_a(r)\to \Omega_p$ is Krasner analytic on
$D_a(r)$. Let  $(c_k)$ in $\Omega_p$ be such that $\lim\limits_{k\to\infty}r^k|c_k|_p=0$ and  all $x\in D_a(r)$,
$$f(x)=\sum\limits^\infty_{k=0}c_k(x-a)^k.$$
We define $||f||_r=\max\limits_k\,r^k|c_k|_p$, then
$\max\limits_{x\in D_a(r)}|f(x)|_p$ is attained when $|x-a|_p=r$ and equals $||f||_r.$

\noindent
{\bf Proof.} See {\bf [K]}: Lemma 3, p. 130.

\medskip
\proclaim 1.12 Lemma.  Let $\emptyset\not=\sigma\subseteq\Omega_p$ be compact. Let $0<r\in |\Omega_p|_p$ and let
$\Gamma\in\Omega_p$, with $|\Gamma|_p=r$.
Assume that
$$\eqalignno{D_\sigma(r)&=\bigcup\limits^M_{i=1}D_{b_i}(r); &(1)\cr
           \sigma\subseteq &\bigcup\limits^M_{i=1}D_{b_i}(r^-),\cr}$$
where $b_1,\ldots,b_M\in\sigma$ and the $D_{b_i}(r)$ are disjoint. Let $\varphi\in H_0(\overline{\sigma})$. Suppose that
$a_1,\ldots,a_N\in\Omega_p$, with
$$ \sigma\subseteq \bigcup\limits^N_{i=1}D_{a_i}(r^-), \eqno (2)$$
where the $D_{a_i}(r)$ are disjoint. Finally, assume that $f$ is Krasner analytic on each of the $D_{a_i}(r)$.
Then $f\in B_r(\sigma)$ and
$$\sum\limits^M_{i=1}\,\int\limits_{b_i,\Gamma}\,f(x)(x-b_i)\varphi(x)\,dx =
  \sum\limits^N_{i=1}\,\int\limits_{a_i,\Gamma}\,f(x)(x-a_i)\varphi(x)\,dx.$$

\medskip
\noindent
{\bf Proof.} By (2) and Lemma 1.4, there exists $\emptyset\not=I\subseteq \{\,1,\ldots,N\,\}$ such that
$$ D_\sigma(r^-)=\bigcup\limits_{i\in I}D_{a_i}(r^-),\quad  D_\sigma(r)=\bigcup\limits_{i\in I}D_{a_i}(r).\eqno (3)$$
To show that $f\in B_\sigma(r)$, let $b\in \Omega_p$, with $D_b(r)\subseteq D_\sigma(r)$. Then (3) implies that
there exists $i\in I$ such that $b\in D_{a_i}(r)$, and hence $D_b(r)=D_{a_i}(r)$; because $f$ is Krasner analytic on
$D_{a_i}(r)$, we see that $f$ is Krasner analytic on $D_b(r)$. This proves that $f\in B_\sigma(r)$.
Now, by (1), (3) and Lemma 1.10 we have
$$\sum\limits^M_{i=1}\,\int\limits_{b_i,\Gamma}\,f(x)(x-b_i)\varphi(x)\,dx=
  \sum\limits_{i\in I}\int\limits_{a_i,\Gamma}\,f(x)(x-a_i)\varphi(x)\,dx.\eqno (4)$$
Now let $1\le j\le N$, with $j\not\in I$. We claim that $D_{a_j}(r)\subseteq \overline{\sigma}$. For suppose that
$b\in\sigma$, with $b\in D_{a_j}(r)$. Then because $\sigma\subseteq  D_\sigma(r)$, (3) implies that 
$b\in D_{a_i}(r)$ for some $i\in I$. We have $i\not=j$ and
$b\in D_{a_i}(r)\cap D_{a_j}(r)$. This contradicts the assumption that $D_{a_i}(r)$ and $ D_{a_j}(r)$ are disjoint
when $i\not=j$. This proves the claim. By assumption, $\varphi\in H_0(\overline{\sigma})$, and since $f\in B_\sigma(r)$, 
we see that $f(x)(x-a_j)\varphi(x)$ is Krasner analytic on $D_{a_j}(r)$ for each $1\le j\le N$ such that $j\not\in I$.
It follows from the $p$-adic Cauchy integral formula ({\bf [K]}: Lemma 4, p. 131) that for $j\not\in I$,
$\int\limits_{a_j,\Gamma}\,f(x)(x-a_j)\varphi(x)\,dx=0$. Therefore, (4) gives
$$\eqalign{\sum\limits^M_{i=1}\,\int\limits_{b_i,\Gamma}\,f(x)(x-b_i)\varphi(x)\,dx=&
  \sum\limits_{i\in I}\int\limits_{a_i,\Gamma}\,f(x)(x-a_i)\varphi(x)\,dx\cr
      = &\sum\limits_{i\in I}\int\limits_{a_i,\Gamma}\,f(x)(x-a_i)\varphi(x)\,dx+\cr
                &\sum\limits_{i\not\in I}\int\limits_{a_i,\Gamma}\,f(x)(x-a_i)\varphi(x)\,dx\cr
		               = &\sum\limits^N_{i=1}\,\int\limits_{a_i,\Gamma}\,f(x)(x-a_i)\varphi(x)\,dx.\cr}$$
This completes the proof of the lemma. \vrule height 6pt width 5pt depth 4pt

\medskip
\proclaim 1.13 Definition. A {\bf $p$-adic Banach space} over $\Omega_p$ is 
a vector space $\cal X$ over $\Omega_p$ together with a
norm $\|\cdot\|_p$ from $\cal X$ to the nonnegative real numbers 
such that for all $x,y\in{\cal X}$: {\rm (a)} $\|x\|_p=0$ if and 
only if $x=0$; {\rm (b)} $\|x+y\|_p\le \max\{\,\|x\|_p,\|y\|_p\,\}$; i
{\rm (c)} $\|ax\|_p=|a|_p\|x\|_p$; {\rm (d)} $\cal X$ is complete
under $\|\cdot\|_p$. We shall assume that $\|{\cal X}\|_p=|\Omega_p|_p$, i.e., for every $x\not=0$ in $\cal X$ there 
exists $a\in\Omega_p$ such that $\|ax\|_p=1$. The {\bf dual} ${\cal X}^*$ of a  $p$-adic Banach space 
over $\Omega_p$ is defined in the usual way. A  {\bf $p$-adic Banach algebra} over  $\Omega_p$ is a $p$-adic Banach space  
$\cal A$ over $\Omega_p$ such that for $x,y\in{\cal A}$, we have $\|x y\|_p\le \|x\|_p \|y\|_p$. We shall assume that
$\cal A$ has a unit. For $x\in{\cal A}$, the {\bf spectrum} $\sigma_x$ of $x$ has the usual meaning, and the 
{\bf resolvent} of $x$ is defined by $R(z;x)=(z-x)^{-1}$, $z\not\in\sigma_x$. If $\cal X$ and $\cal Y$ are $p$-adic
Banach spaces over $\Omega_p$, then $B({\cal X},{\cal Y})$ is the vector space of $\Omega_p$-linear continuous maps
form $\cal X$ to $\cal Y$. $B({\cal X},{\cal Y})$ is a $p$-adic Banach space under the usual operator norm, and 
$B({\cal X})=B({\cal X},{\cal X})$ is a $p$-adic Banach algebra under this operator norm. 

\medskip
\proclaim 1.14 Definition. Let $\cal X$ be a $p$-adic Banach space over $\Omega_p$. An operator $A\in B({\cal X})$, with
compact spectrum $\sigma_A$, is called {\bf analytic} if the resolvent $R(z;A)$ is {\bf Krasner analytic} on
$\overline{\sigma}_A$, in the sense that for all $x\in{\cal X}$ and $h\in{\cal X}^*$, 
the function $z\mapsto h(R(z;A)x)$ is in $H_0(\overline{\sigma}_A)$.

\medskip
\proclaim 1.15 Definition. Let $J$ be any nonempty indexing set, and define $\Omega_p(J)$ to be the set of all
sequences $c=\displaystyle{(c_j)_{j\in J}}$ in $\Omega_p$ such that for every $\epsilon>0$ only finitely many
$|c_j|_p$ are $>\epsilon$. Define $||c||_p=\max\limits_j|c_j|_p$. Then $\Omega_p(J)$ is a $p$-adic Banach space over
$\Omega_p$. 

\medskip
\proclaim 1.16 Lemma.  Let ${\cal X}\simeq\Omega_p(J)$, with $J\not=\emptyset$. Let $\sigma\subseteq\Omega_p$
be compact. Let $F:\overline\sigma\to B({\cal X})$ be an analytic operator-valued function, i.e., for all
$y\in{\cal X},h\in {\cal X}^*$ the $\Omega_p$-valued function
$$F_{h,y}(x)=h(F(x)y),\quad x\in\overline\sigma$$
belongs to $H_0(\overline\sigma)$. Let $a\in\Omega_p$, and $0<r\in|\Omega_p|_p$. Let $\Gamma\in\Omega_p$, with
$|\Gamma|_p=r$. 
Assume that there is no $b\in\sigma$ such that $|b-a|_p=r$. Finally, let
$f:D_a(r)\to\Omega_p$ be Krasner analytic on $D_a(r)$. Define
$$S_n={1\over n}\sum\limits_{\xi^n=1}f(a+\xi\Gamma)F(a+\xi\Gamma),\quad n=1,2,\ldots.$$
Then the limit
$$\lim\limits_{\buildrel{n\to\infty}\over{p\not\;\mid n}}S_n$$
exists as an operator in $B({\cal X})$. We define $\int\limits_{a,\Gamma}\,f(x)F(x)\,dx$ to be this operator.
Moreover, for any $y\in{\cal X}$ and $h\in{\cal X}^*$, we have
$$h\bigg(\bigg(\int\limits_{a,\Gamma}\,f(x)F(x)\,dx\bigg)(y)\bigg)=\int\limits_{a,\Gamma}\,f(x)F_{h,y}(x)\,dx.$$
In particular, assume that an analytic operator $A\in B({\cal X})$ has compact spectrum $\sigma_A$.
Then by definition, the resolvent $R(x;A)=(x-A)^{-1}$ is Krasner analytic on $\overline{\sigma}_A$
and hence the integral
$$\int\limits_{a,\Gamma}\,f(x)(x-a)R(x;A)\,dx=\int\limits_{a,\Gamma}\,f(x)(x-a)R(x;A)\,dx$$
exists, such that for all $y\in{\cal X}$ and all $h\in{\cal X}^*$, if  $R_{h,y}(x;A)=h(R(x;A)y)$, then
$$h\bigg(\bigg(\int\limits_{a,\Gamma}\,f(x)(x-a)R(x;A)\,dx\bigg)(y)\bigg)
            =h\bigg(\bigg(\int\limits_{a,\Gamma}\,f(x)(x-a)R_{h,y}(x;A)\,dx\bigg)(y)\bigg).$$

\medskip
\noindent
{\bf Proof.} See {\bf [K]}: Lemma 3, p. 149. 

\medskip
\proclaim 1.17 Definition. Let $\cal A$ be a $p$-adic Banach algebra over $\Omega_p$. Let $A\in{\cal A}$ have
spectrum $\sigma_A\not=\emptyset$. Define ${\cal F}(A)=L(\sigma_A)$. 

\medskip 
\proclaim 1.18 Spectral Theorem I (Vishik). Let ${\cal X}\simeq\Omega_p(J)$, where $J$ is a nonempty indexing set. Let $A\in B({\cal X})$ have compact
spectrum $\sigma_A\not=\emptyset$, and assume that $A$ is analytic. Let $f\in{\cal F}(A)$. 
Let $0<r\in |\Omega_p|_p$ be such 
that $f\in B_r(\sigma_A)$, where
$$D_{{\sigma}_A}(r)=\bigcup\limits^N_{i=1}D_{a_i}(r),\quad \sigma_A\subseteq D_{{\sigma}_A}(r^-).$$
Here the union is assumed to be disjoint and each $a_i\in\sigma_A$. Let $\Gamma\in\Omega_p$, with $|\Gamma|_p=r$. 
Define 
$$f(A)=\sum\limits^N_{i=1}\int\limits_{a_i,\Gamma}\,f(x)(x-a_i)R(x;A)\,dx.$$ 
Let $g\in {\cal F}(A)$, and let 
$\alpha,\beta\in\Omega_p$. Then \item\item{(1)} $(\alpha f+\beta g)(A)=\alpha f(A)+\beta g(A)$;
\item\item{(2)} $(f\cdot g)(A)=f(A)g(A)$;
\item\item{(3)} The mapping $f\mapsto f(A)$ is continuous on $B_r(\sigma_A)$.

\medskip
\noindent
{\bf Proof.} See {\bf [K]}: Spectral Theorem, p. 149.

\medskip
\proclaim 1.19 Definition.
Let $\cal A$ be a $p$-adic Banach algebra and let $\cal X$ be a $p$-adic Banach space.
Let $A\in{\cal A}$ have compact spectrum $\sigma_A$. Then we say that $A$ is ${\cal X}$-{\bf analytic} 
iff there exists a continuous $\Omega_p$-monomorphism $\theta:{\cal A}\to B({\cal X})$ such that
$\theta^{-1}:\theta({\cal A})\to {\cal A}$ is continuous and $\theta(A)$ is analytic in $B({\cal X})$. Such a 
monomorphism  $\theta:{\cal A}\to B({\cal X})$ is said to be an {\bf embedding} of $\cal A$ into $B({\cal X})$.
Let $B\in{\cal A}$. Then $B$ is $\cal A$-{\bf analytic} iff there exists a sequence $(B_k)$ in $\cal A$ and a sequence
$({\cal X}_k)$ of $p$-adic Banach spaces ${\cal X}_k$ of the form 
${\cal X}_k\simeq\Omega_p(J_k)$, with $J_k\not=\emptyset$, such that $\sigma_{B_k}\not=\emptyset$ is compact, 
$B_k$ is ${\cal X}_k$-analytic and $\lim\limits_{k\to\infty}\|B-B_k\|_p=0$.

\medskip
\proclaim 1.20 Lemma. Let $\cal A$ and $\cal X$ be as in Definition 1.19, with  ${\cal X}\simeq\Omega_p(J)$, 
where $J$ is a nonempty indexing set. Let $A\in{\cal A}$ have compact spectrum 
$\emptyset\not=\sigma_A$, and assume that $A$ is $\cal X$-analytic. Let $\theta:{\cal A}\to B({\cal X})$ be an embedding
of $\cal A$ into $B({\cal X})$. Let $0<r\in |\Omega_p|_p$ and
$\Gamma\in\Omega_p$, with $|\Gamma|_p=r$.
Assume that $f\in B_{\sigma_A}(r)$ and
$$\sigma\subseteq \bigcup\limits^N_{i=1}D_{a_i}(r^-),$$ 
where the $D_{a_i}(r)$ are disjoint. Then we have
$$\eqalignno{\theta^{-1}[f(\theta(A))]&=
\lim\limits_{\buildrel{n\to\infty}\over{p\not\;\mid n}}
\sum\limits^N_{i=1}{1\over n}\sum\limits_{\xi^n=1}f(x^i_\xi)(x^i_\xi-a_i)R(x^i_\xi;A),&(1)\cr 
                                    &=\sum\limits^N_{i=1}\int\limits_{a_i,\Gamma}f(x)(x-a_i)R(x;A)\,dx,}$$
where for $1\le i\le N$ and $\xi\in\Omega_p$, we define $x^i_\xi=a_i+\xi\Gamma$. 
Moreover, the limit in (1) is independent of the $\Gamma,r,a_1,\ldots,a_N$. Hence we may
define $f(A)$ by $f(A)=\theta^{-1}[f(\theta(A))]$, and this definition will not depend on $\theta$. 

\medskip
\noindent
{\bf Proof.} We have $\sigma_{\theta(A)}=\sigma_A$, hence by Theorem 1.18 (Vishik's Spectral Theorem I), we have
$$\eqalign{\theta^{-1}[f(\theta(A))]&=
                          \theta^{-1}\bigg[\sum\limits^N_{i=1}\int\limits_{a_i,\Gamma}f(x)(x-a_i)R(x;A)\,dx\bigg]\cr
                      &=\theta^{-1}\bigg[\lim\limits_{\buildrel{n\to\infty}\over{p\not\;\mid n}}
\sum\limits^N_{i=1}{1\over n}\sum\limits_{\xi^n=1}f(x^i_\xi)(x^i_\xi-a_i)R(x^i_\xi;\theta(A))\bigg]\cr
&=\lim\limits_{\buildrel{n\to\infty}\over{p\not\;\mid n}} 
\sum\limits^N_{i=1}{1\over n}\sum\limits_{\xi^n=1}f(x^i_\xi)(x^i_\xi-a_i)\theta^{-1}[R(x^i_\xi;\theta(A))]\cr
&=\lim\limits_{\buildrel{n\to\infty}\over{p\not\;\mid n}}
\sum\limits^N_{i=1}{1\over n}\sum\limits_{\xi^n=1}f(x^i_\xi)(x^i_\xi-a_i)R(x^i_\xi;A).\cr}$$
Finally, by Vishik's Spectral Theorem I, the representation
$$f(\theta(A))=\sum\limits^N_{i=1}\int\limits_{a_i,\Gamma}\,f(x)(x-a_i)R(x;\theta(A))\,dx$$
does not depend on $\Gamma,r,a_1,\ldots,a_N$. \vrule height 6pt width 5pt depth 4pt

\medskip
\proclaim 1.21 Spectral Theorem II (Vishik). Let ${\cal X}\simeq\Omega_p(J)$, where $J$ is a nonempty indexing set. Let 
Let $\cal A$ be a $p$-adic Banach algebra, and let $A\in {\cal A}$ have compact spectrum $\sigma_A\not=\emptyset$,
and assume that $A$ is $\cal X$-analytic. Let $f\in{\cal F}(A)$. Let $0<r\in |\Omega_p|_p$ be such
that $f\in B_r(\sigma_A)$, where
$$D_{{\sigma}_A}(r)=\bigcup\limits^N_{i=1}D_{a_i}(r),\quad \sigma_A\subseteq D_{{\sigma}_A}(r^-).$$
Here the union is assumed to be disjoint and each $a_i\in\sigma_A$. Let $\Gamma\in\Omega_p$, with $|\Gamma|_p=r$.
Define
$$f(A)=\sum\limits^N_{i=1}\int\limits_{a_i,\Gamma}\,f(x)(x-a_i)R(x;A)\,dx.$$
Let $g\in {\cal F}(A)$, and let
$\alpha,\beta\in\Omega_p$. Then \item\item{(1)} $(\alpha f+\beta g)(A)=\alpha f(A)+\beta g(A)$;
\item\item{(2)} $(f\cdot g)(A)=f(A)g(A)$;
\item\item{(3)} The mapping $f\mapsto f(A)$ is continuous $B_r(\sigma_A)$.

\medskip
\noindent
{\bf Proof.} Let $\theta:{\cal A}\to B({\cal X})$ be an embedding of $\cal A$ into $B({\cal X})$, where
${\cal X}\simeq\Omega_p(J)$, $J\not=\emptyset$. Assume that $\theta(A)$ is $\cal X$-analytic. Then Theorem 1.18 applies
to $\theta(A)$, and hence Theorem 1.21 follows form Lemma 1.20 and the fact that for $f\in B_{\sigma_A}(r)$, 
$f(A)=\theta^{-1}[f(\theta(A))]$. \vrule height 6pt width 5pt depth 4pt

\vskip 1cm
\centerline{\bf 2. EXTENDED $p$-ADIC SPECTRAL THEORY}
\vskip 1cm
\proclaim 2.1 Theorem (Perturbation Theory). Let ${\cal A}$ be a $p$-adic Banach algebra. Let
$A\in{\cal A}$ spectrum $\sigma_A$. Let $\epsilon>0$ be any positive real number. 
Then there exists a constant $N_{\epsilon}>0$ such that
$$\|R(x;A)\|_p\le N_{\epsilon}$$
for all $x\not\in D_{\sigma_A}(\epsilon^-)$. Moreover, there exists a $\delta>0$ such that if
$B\in {\cal A}$ and $||A-B||_p<\delta$, then
$$\sigma_B\subseteq D_{\sigma_A}(\epsilon^-),$$
and $$||R(x;A)-R(x;B)||_p<\epsilon  \hbox{ for all } x\not\in D_{\sigma_A}(\epsilon^-).$$

\noindent
{\bf Proof.} See {\bf [DS]}: Lemma 3, p. 585. The proof presented in {\bf [DS]} 
works for $p$-adic Banach algebras, and the
proof does not assume that $\sigma_A$ is compact, nor does it assume that $\sigma_A\not=\emptyset$. 
\vrule height 6pt width 5pt depth 4pt

\medskip
\proclaim Remark.  Let $\cal A$ be a $p$-adic Banach algebra, and let $A\in{\cal A}$ have spectrum $\sigma_A$.
If $A$ is analytic, then $\sigma_A\not=\emptyset$. To see this, let $\epsilon>0$ be arbitrary, then by Lemma 2.1,
there exists a $\delta>0$ such that if $B\in{\cal A}$ and $\|A-B\|_p<\delta$, then
$$\sigma_B\subseteq D_{\sigma_A}(\epsilon^-). \eqno (*)$$
Because $\lim\limits_{k\to\infty}\|A-A_k\|_p=0$, there exists a $k$ such that $\|A-A_k\|_p<\delta$. Then by (*),
we have
$$\emptyset\not=\sigma_{A_k}\subseteq D_{\sigma_A}(\epsilon^-).$$
Hence, $D_{\sigma_A}(\epsilon^-)\not=\emptyset$, which implies that $\sigma_A\not=\emptyset$.

\medskip
\proclaim 2.2 Lemma. Let $\cal A$ be a $p$-adic Banach algebra, and let $A\in{\cal A}$ have compact spectrum
$\sigma_A\not=\emptyset$. Let $\cal X$ be a $p$-Banach space such that ${\cal X}\simeq\Omega_p(J)$, with 
$J\not=\emptyset$, and suppose that $A$ be $\cal X$-analytic with spectrum $\sigma_A$. 
Let $0<r_1\le r_2$, and let $\Gamma_1,\Gamma_2$ be in 
$\Omega_p$,  with $|\Gamma_1|_p=r_1$, $|\Gamma_2|_p=r_2$. Assume that
$$\eqalignno{D_\sigma(r_1)&=\bigcup\limits^M_{i=1}D_{a_i}(r_1),\quad \sigma \subseteq
                          \bigcup\limits^M_{i=1}D_{a_i}({r_1}^-);& (1)\cr
           D_\sigma(r_2)&=\bigcup\limits^N_{i=1}D_{b_i}(r_2), \quad \sigma \subseteq
                          \bigcup\limits^N_{i=1}D_{b_i}({r_2}^-),\cr}$$
where the $D_{a_i}(r_1),D_{b_i}(r_2)$ are disjoint.
If $f\in B_{r_2}(\sigma_A)$, then $f\in B_{r_1}(\sigma_A)$, and the following
sums exist and are equal.
$$\sum\limits^M_{i=1}\,\int\limits_{a_i,\Gamma_1}\,f(x)(x-a_i)R(x;A)\,dx=
  \sum\limits^N_{i=1}\,\int\limits_{b_i,\Gamma_2}\,f(x)(x-b_i)R(x;A)\,dx.\eqno (2)$$

\medskip
\noindent
{\bf Proof.} Let $\theta:{\cal A}\to B({\cal X})$ be and embedding of ${\cal A}$ into $B({\cal X})$ such that 
$\theta(A)$ is analytic. Let $y\in{\cal X}$ and $h\in{\cal X}^*$. Then the function 
$\varphi_{h,y}(x)=h(R(x;\theta(A))y)$ is in $H_0(\overline{\sigma}_A)$. Let $f\in B_{r_2}(\sigma_A)$, then it is clear
that  $f\in B_{r_1}(\sigma_A)$. Then by Lemma 1.10, (1) implies that 
$$\sum\limits^M_{i=1}\,\int\limits_{a_i,\Gamma_1}\,f(x)(x-a_i))\varphi_{h,y}(x)\,dx=
  \sum\limits^N_{i=1}\,\int\limits_{b_i,\Gamma_2}\,(x)f(x)(x-b_i)\varphi_{h,y}(x)\,dx.\eqno (3)$$
Because $y\in{\cal X}$ and $h\in{\cal X}^*$ are arbitrary, by Lemma 1.16, (3) implies that
$$\sum\limits^M_{i=1}\,\int\limits_{a_i,\Gamma_1}\,f(x)(x-a_i)R(x;\theta(A))\,dx=
  \sum\limits^N_{i=1}\,\int\limits_{b_i,\Gamma_2}\,f(x)(x-b_i)R(x;\theta(A))\,dx.\eqno (4)$$
Because $\theta$ is bicontinuous, (4) implies (2). \vrule height 6pt width 5pt depth 4pt

\medskip
\proclaim 2.3 Lemma. Let $\cal A$ be a $p$-adic Banach algebra, and let $A\in{\cal A}$ have compact spectrum
$\sigma_A\not=\emptyset$. Let $\cal X$ be a $p$-Banach space such that ${\cal X}\simeq\Omega_p(J)$, with
$J\not=\emptyset$, and suppose that $A$ be $\cal X$-analytic. Let $r>0$ and let $\Gamma\in\Omega_p$, with $|\Gamma|_p=r$.
Assume that 
$$D_{\sigma_A}(r)=\bigcup\limits^M_{i=1}D_{b_i}(r),\quad
              \sigma_A\subseteq \bigcup\limits^M_{i=1}D_{b_i}(r^-),$$
where $b_1,\ldots,b_M$ are in $\Omega_p$ and the $D_{b_i}(r)$ are disjoint. Finally, suppose that
$$\sigma_A\subseteq \bigcup\limits^N_{i=1}D_{a_i}(r^-),$$
where $a_1,\ldots,a_N$ are in $\Omega_p$ and the $D_{a_i}(r)$ are disjoint. 
If $f$ is Krasner analytic on each $D_{a_i}(r)$, then we have $f\in B_r(\sigma_A)$ and
$$\sum\limits^M_{i=1}\,\int\limits_{b_i,\Gamma}\,f(x)(x-b_i)R(x;A)\,dx =
\sum\limits^N_{i=1}\,\int\limits_{a_i,\Gamma}\,f(x)(x-a_i)R(x;A)\,dx.$$

\noindent
{\bf Proof.} By Lemma 1.12, the proof of Lemma 2.3 is similar to the proof of \hbox{Lemma 2.2.
\vrule height 6pt width 5pt depth 4pt}

\medskip
\proclaim 2.4 Lemma. Let $\cal A$ be a $p$-adic Banach algebra. Let $A\in{\cal A}$ be $\cal A$-analytic.
Let $(A_k)$ be a sequence in $\cal A$ and $({\cal X}_k)$ a sequence of $p$-adic Banach spaces ${\cal X}_k$ of the form
${\cal X}_k\simeq\Omega_p(J_k)$, with $J_k\not=\emptyset$, such that $\sigma_{A_k}\not=\emptyset$ is compact,
$A_k$ is ${\cal X}_k$-analytic and $\lim\limits_{k\to\infty}\|A-A_k\|_p=0$ (see Definition 1.19). Let $A$ have
spectrum $\sigma_A$. Let $r>0$, with $r\in|\Omega_p|_p$, and let $\Gamma\in \Omega_p$, with $|\Gamma|_p=r$.
Assume that
$$ \sigma_A\subseteq \bigcup\limits^N_{i=1}D_{a_i}(r^-),\eqno (1)$$
where $a_1,\ldots,a_N\in\Omega_p$ and the $D_{a_i}(r)$ are disjoint. Assume that for each $k$ there exists 
$a_{k,1},\ldots,a_{k,N_k}$ in $\sigma_{A_k}$ such that 
$$\eqalignno{D_{\sigma_{A_k}}(r)&=\bigcup\limits^{N_k}_{i=1}D_{a_{k,i}}(r); &(2)\cr
                \sigma_{A_k}&\subseteq \bigcup\limits^{N_k}_{i=1}D_{a_{k,i}}(r^-),\cr}$$
where the $D_{a_{k,i}}(r)$ are disjoint. Then $\sigma_A\not=\emptyset$. Moreover, if $f$ is Krasner analytic on the
$D_{a_i}(r)$,then the following limits exist and are equal.
$$\lim\limits_{k\to\infty}\sum\limits^N_{i=1}\,\int\limits_{a_i,\Gamma}\,f(x)(x-a_i)R(x;A_k)\,dx=
   \lim\limits_{k\to\infty}\sum\limits^{N_k}_{i=1}\,\int\limits_{a_{k,i},\Gamma}\,f(x)(x-a_{k,i})R(x;A_k)\,dx.$$

\medskip
\noindent
{\bf Proof.} 
Because $f$ is Krasner analytic on each $D_{a_i}(r)$,
Lemma 1.11 implies that for each $1\le i\le N$, $\max\limits_{x\in D_{a_i}(r)}|f(x)|_p$ 
is attained when $|x-a_i|_p=r$, and hence we may write
$$\rho_i=\max\limits_{x\in D_{a_i}(r)}|f(x)|_p=\max\limits_{|x-a_i|_p=r}|f(x)|_p.$$
Now let $0<\epsilon<r$ be arbitrary. By Lemma 2.1, there exists a positive number $\delta>0$ such that if
$B\in {\cal A}$ and $\|A-B\|_p<\delta$, then
$$\sigma_B\subseteq D_{\sigma_A}(\epsilon^-),\eqno (3)$$
and
$$\|R(x;A)-R(x;B)\|_p<{\epsilon\over r(1+\max\limits_j\rho_j)},\quad x\not\in D_{\sigma_A}(\epsilon^-).$$
Now let $k_0$ be so large that $\|A-A_k\|_p<\delta$ for $k\ge k_0$. Then for $k\ge k_0$, we have
we have $\sigma_{A_k}\subseteq D_{\sigma_A}(\epsilon^-)\subseteq D_{\sigma_A}(r^-)$. 
By Lemma 1.3, (1) implies that there exists 
$\emptyset\not=I\subseteq \{\,1,\ldots,N\,\}$ such that
$$D_{\sigma_A}(r^-)=\bigcup\limits_{i\in I}D_{a_i}(r^-)\subseteq \bigcup\limits^N_{i=1}D_{a_i}(r^-). \eqno (4)$$
For $k\ge k_0$, (4) gives
$$\sigma_{A_k}\subseteq D_{\sigma_A}(r^-)\subseteq \bigcup\limits^N_{i=1}D_{a_i}(r^-),$$
and hence by Lemma 2.3, (1) and (2) imply that
$$\sum\limits^{N_k}_{i=1}\int\limits_{a_{k,i},\Gamma}f(x)(x-a_{k,i})R(x;A_k)\,dx=
  \sum\limits^N_{i=1}\int\limits_{a_i,\Gamma}f(x)(x-a_i)R(x;A_k)\,dx.$$
Therefore if $k\ge k_0$ we get, with $x^i_{\xi}=a_i+\xi\Gamma$, 
$$\eqalignno{&\sum\limits^{N_k}_{i=1}\int\limits_{a_{k,i},\Gamma}f(x)(x-a_{k,i})R(x;A_k)\,dx-
           \sum\limits^{N_l}_{i=1}\int\limits_{a_{l,i},\Gamma}f(x)(x-a_{l,i})R(x;A_l)\,dx= &(5)\cr
	   &\sum\limits^N_{i=1}\int\limits_{a_i,\Gamma}f(x)(x-a_i)R(x;A_k)\,dx-
	     \sum\limits^N_{i=1}\int\limits_{a_i,\Gamma}f(x)(x-a_i)R(x;A_l)\,dx=\cr
           \lim\limits_{\buildrel{n\to\infty}\over{p\not\;\mid n}}
&\sum\limits^N_{i=1}{1\over n}\sum\limits_{\xi^n=1}
\big[f(x^i_{\xi})(x^i_{\xi}-a_i)(R(x^i_{\xi};A_k)-R(x^i_{\xi};A))\big]-\cr
\lim\limits_{\buildrel{n\to\infty}\over{p\not\;\mid n}}&\sum\limits^N_{i=1}{1\over n}\sum\limits_{\xi^n=1}
\big[f(x^i_{\xi})(x^i_{\xi}-a_i)(R(x^i_{\xi};A_l)-R(x^i_{\xi};A))\big].}$$
Now, we claim that for $1\le i\le N$ and $\xi\in\Omega_p$, with, $|\xi|_p=1$,
$$x^i_{\xi}\not\in D_{\sigma_A}(\epsilon^-).\eqno (6)$$
To prove this claim, suppose that $x^i_{\xi}\in D_{\sigma_A}(\epsilon^-)$, then  $x^i_{\xi}\in D_{\sigma_A}(r^-)$, hence
(4) implies that for some $1\le j\le N$, $x^i_{\xi}\in D_{a_j}(r^-)$. We must have $i=j$, for if $i\not=j$, then
because $|x^i_{\xi}-a_i|_p=|\xi \Gamma|_p=r$, we would have $x^i_{\xi}\in D_{a_i}(r)\cap D_{a_j}(r)=\emptyset$. Thus,
$i=j$, which implies that $x^i_{\xi}\in D_{a_i}(r^-)$. But then we have the contradiction that
$r=|x^i_{\xi}-a_i|_p<r$. Thus (6) holds.
If $p\not\;\mid n$, $1\le i\le N$ and $m\ge k_0$, then (3) and (6) imply that 
$$\eqalignno{& \bigg\|{1\over n}\sum\limits_{\xi^n=1}
f(x^i_{\xi})(x^i_{\xi}-a_i)(R(x^i_{\xi};A_m)-R(x^i_{\xi};A))\bigg\|_p &  (7)\cr
 \le & \max\limits_{|x-a_i|_p=r}|f(x)|_p \cdot |\Gamma|_p \cdot \max\limits_{\xi^n=1}\|R(x^i_\xi;A_m)-R(x^i_\xi;A)\|_p\cr
	                \le &(r \rho_i)\bigg({\epsilon\over r(1+\max\limits_j\rho_j)}\bigg)\cr
			            <&\epsilon.\cr}$$
Here we have used the fact that $D_{\sigma_A}(\epsilon^-)\subseteq D_{\sigma_A}(r^-)$, and hence
$x\not\in D_{\sigma_A}(r^-)$ implies $x\not\in D_{\sigma_A}(\epsilon^-)$.
Therefore, (5) and (7) imply that for $k,l\ge k_0$, 
$$\bigg\|\sum\limits^{N_k}_{i=1}\int\limits_{a_{k,i},\Gamma}f(x)(x-a_{k,i})R(x;A_k)\,dx-
           \sum\limits^{N_l}_{i=1}\int\limits_{a_{l,i},\Gamma}f(x)(x-a_{l,i})R(x;A_l)\,dx\bigg\|_p
            \le\epsilon.$$
Because $0<\epsilon<r$ is arbitrary, we see that for $k\ge k_0$ the terms
$$\sum\limits^{N_k}_{i=1}\int\limits_{a_{k,i},\Gamma}f(x)(x-a_{k,i})R(x;A_k)\,dx=
  \sum\limits^N_{i=1}\int\limits_{a_i,\Gamma}f(x)(x-a_i)R(x;A_k)\,dx$$
form a Cauchy sequence. Hence the following limits exist and are equal.
$$\lim\limits_{k\to\infty}\sum\limits^{N_k}_{i=1}\int\limits_{a_{k,i},\Gamma}f(x)(x-a_{k,i})R(x;A_k)\,dx=
  \lim\limits_{k\to\infty}\sum\limits^N_{i=1}\int\limits_{a_i,\Gamma}f(x)(x-a_i)R(x;A_k)\,dx.$$
This completes the proof of the lemma. \vrule height 6pt width 5pt depth 4pt

\medskip
\proclaim 2.5 Lemma. Let $\cal A$ be a $p$-adic Banach algebra. Let $A\in{\cal A}$ be $\cal A$-analytic. Let
$(A_k)$ be a sequence in $\cal A$ and $({\cal X}_k)$ a sequence of $p$-adic Banach spaces ${\cal X}_k$ of the form
${\cal X}_k\simeq\Omega_p(J_k)$, with $J_k\not=\emptyset$, such that $\sigma_{A_k}\not=\emptyset$ is compact,
$A_k$ is ${\cal X}_k$-analytic and $\lim\limits_{k\to\infty}\|B-B_k\|_p=0$. Let $A$ have
spectrum $\sigma_A$. Let $0<r_2\le r_1$ be in $|\Omega_p|_p$, and let $\Gamma_1,\Gamma_2$ be in $\Omega_p$, with 
$|\Gamma_1|_p=r_1$, $|\Gamma_2|_p=r_2$. Assume that $a_1,\ldots,a_M;b_1,\ldots,b_N$ in $\Omega_p$ are given, with
$$\sigma_A\subseteq\bigcup\limits^M_{i=1}D_{a_i}(r^-_1)\hbox{ and }
\sigma_A\subseteq\bigcup\limits^N_{i=1}D_{b_i}({r_2}^-), \eqno (1)$$
where the $D_{a_i}(r_1)$ are disjoint and the $D_{b_i}(r_2)$ are disjoint. Let $f$ be Krasner analytic on the
$D_{a_i}(r_1)$ and on the  $D_{b_i}(r_2)$. Then the following limits exist and are equal:
$$\lim\limits_{k\to\infty}\sum\limits^M_{i=1}\,\int\limits_{a_i,\Gamma_1}\,f(x)(x-a_i)R(x;A_k)\,dx=
\lim\limits_{k\to\infty}\sum\limits^N_{i=1}\,\int\limits_{b_i,\Gamma_2}\,f(x)(x-b_i)R(x;A_k)\,dx.$$

\noindent
{\bf Proof.} By (1) and Lemma 2.1, we may assume that for all $k$,
$$\sigma_{A_k}\subseteq \bigcup\limits^M_{i=1}D_{a_i}(r^-_1) \hbox{ and } 
\sigma_{A_k}\subseteq \bigcup\limits^N_{i=1}D_{b_i}(r^-_2).\eqno (2)$$
Then by Lemma 1.3, (2) implies that for all $k$ there exist $\emptyset\not= I_k,J_k$ such that
$I_k\subseteq \{\,1,\ldots,M\,\}$, $J_k\subseteq \{\,1,\ldots,N\,\}$
and
$$\eqalignno{\sigma_{A_k}\subseteq D_{\sigma_{A_k}} (r_1^-)&=\bigcup\limits_{i\in I_k}D_{a_i}(r_1^-)\hbox{ and }
D_{\sigma_{A_k}} (r_1)=\bigcup\limits_{i\in I_k}D_{a_i}(r_1); & (3)\cr
           \sigma_{A_k}\subseteq D_{\sigma_{A_k}} (r_2^-)&=\bigcup\limits_{i\in J_k}D_{b_i}(r_2^-)\hbox{ and }
	              D_{\sigma_{A_k}} (r_2)=\bigcup\limits_{i\in J_k}D_{b_i}(r_2).\cr}$$
Now let $k$ be arbitrary. By Lemma 2.3, (2) and (3) together imply that
$$\eqalignno{\sum\limits^M_{i=1}\,\int\limits_{a_i,\Gamma_1}\,f(x)(x-a_i)R(x;A_k)\,dx &=
             \sum\limits_{i\in I_k}\;\int\limits_{a_i,\Gamma_1}\,f(x)(x-a_i)R(x;A_k)\,dx; & (4)\cr
             \sum\limits^N_{i=1}\,\int\limits_{b_i,\Gamma_2}\,f(x)(x-b_i)R(x;A_k)\,dx &=
             \sum\limits_{i\in J_k}\;\int\limits_{b_i,\Gamma_2}\,f(x)(x-b_i)R(x;A_k)\,dx.\cr}$$
By Lemma 2.2, (2) and (3) together also imply that
$$\sum\limits_{i\in I_k}\;\int\limits_{a_i,\Gamma_1}\,f(x)(x-a_i)R(x;A_k)\,dx=
  \sum\limits_{i\in J_k}\;\int\limits_{b_i,\Gamma_2}\,f(x)(x-b_i)R(x;A_k)\,dx. \eqno (5)$$
Consequently, (4) and (5) together imply that for all $k$,
$$\sum\limits^M_{i=1}\,\int\limits_{a_i,\Gamma_1}\,f(x)(x-a_i)R(x;A_k)\,dx=
\sum\limits^N_{i=1}\,\int\limits_{b_i,\Gamma_2}\,f(x)(x-b_i)R(x;A_k)\,dx. \eqno (6)$$
It follows from (1), (3), (6), and Lemma 2.4 that $\sigma_A\not=\emptyset$ 
and that the following limits exist and are equal.
$$\lim\limits_{k\to\infty}\sum\limits^M_{i=1}\,\int\limits_{a_i,\Gamma_1}\,f(x)(x-a_i)R(x;A_k)\,dx=
\lim\limits_{k\to\infty}\sum\limits^N_{i=1}\,\int\limits_{b_i,\Gamma_2}\,f(x)(x-b_i)R(x;A_k)\,dx.$$
This completes the proof of the lemma. \vrule height 6pt width 5pt depth 4pt

\medskip
\proclaim 2.6 Lemma. Let $\cal A$ be a $p$-adic Banach algebra. Let $A\in{\cal A}$ be $\cal A$-analytic with spectrum
$\sigma_A$. Let
$(A_k)$ be a sequence in $\cal A$ and $({\cal X}_k)$ a sequence of $p$-adic Banach spaces ${\cal X}_k$ of the form
${\cal X}_k\simeq\Omega_p(J_k)$, with $J_k\not=\emptyset$, such that $\sigma_{A_k}\not=\emptyset$ is compact,
$A_k$ is ${\cal X}_k$-analytic and $\lim\limits_{k\to\infty}\|A-A_k\|_p=0$. Let $A$ have spectrum $\sigma_A$.
Let $0<r\in |\Omega_p|_p$, and let $\Gamma\in\Omega_p$, with $|\Gamma|_p=r$. Assume that $a_1,\ldots,a_N$ $\Omega_p$
are given, with
$$\sigma_A\subseteq\bigcup\limits^N_{i=1}D_{a_i}(r^-), \eqno (1)$$
where the $D_{a_i}(r)$ are disjoint. Assume that $f$ is Krasner analytic on the $D_{a_i}(r)$. Then the following limits
exist and are equal.
$$\lim\limits_{\buildrel{n\to\infty}\over{p\not\;\mid n}}\sum\limits^N_{i=1}{1\over n}
\sum\limits_{\xi^n=1}f(x^i_{\xi})(x^i_{\xi}-a_i)R(x^i_{\xi};A)=
\lim\limits_{k\to\infty}\sum\limits^N_{i=1}\,\int\limits_{a_i,\Gamma}\,f(x)(x-a_i)R(x;A_k)\,dx,$$
where for each $1\le i\le N$ and $\xi\in\Omega_p$, we set $x^i_\xi=a_i+\xi\Gamma$. Therefore the sum
$$\sum\limits^N_{i=1}\,\int\limits_{a_i,\Gamma}\,f(x)(x-a_i)R(x;A)\,dx$$
exits and we have
$$\sum\limits^N_{i=1}\,\int\limits_{a_i,\Gamma}\,f(x)(x-a_i)R(x;A)\,dx=
\lim\limits_{k\to\infty}\sum\limits^N_{i=1}\,\int\limits_{a_i,\Gamma}\,f(x)(x-a_i)R(x;A_k)\,dx.$$

\medskip
\noindent
{\bf Proof.} Because $f$ is Krasner analytic on each of the $D_{a_i}(r)$, it follows from 
Lemma 1.11, for $1\le i\le N$, $\max\limits_{x\in D_{a_i}(r)}|f(x)|_p$ is attained when $|x-a_i|_p=r$, and hence we
may write
$$\rho_i=\max\limits_{x\in D_{a_i}(r)}|f(x)|_p=\max\limits_{|x-a_i|_p=r}|f(x)|_p.$$
Now let $0<\epsilon<r$ be arbitrary. By Lemma 2.1, there exists a positive number $\delta>0$ such that if
$B\in {\cal A}$ and $\|A-B\|_p<\delta$, then
$$\sigma_B\subseteq D_{\sigma_A}(\epsilon^-),$$
and
$$\|R(x;A)-R(x;B)\|_p<{\epsilon\over r(1+\max\limits_j\rho_j)},\quad x\not\in D_{\sigma_A}(\epsilon^-).\eqno (2)$$
Then by Lemma 1.3, (1) implies that for all $k$ there exist $\emptyset\not= I_k$ such that
$I_k\subseteq \{\,1,\ldots,N\,\}$ and 
$$\sigma_{A_k}\subseteq D_{\sigma_{A_k}} (r^-)=\bigcup\limits_{i\in I_k}D_{a_i}(r^-)\hbox{ and }
D_{\sigma_{A_k}} (r)=\bigcup\limits_{i\in I_k}D_{a_i}(r). \eqno  (3)$$
By 2.4, (1) and (3) together imply that the following limit exists.
$$\lim\limits_{k\to\infty}\sum^N_{i=1}\int\limits_{a_i,\Gamma}f(x)(x-a_i)R(x;A_k)\,dx.$$
Now let $k_0$ be so large that $\|A-A_{k_0}\|_p<\delta$ and
$$\bigg \|\sum^N_{i=1}\int\limits_{a_i,\Gamma}f(x)(x-a_i)R(x;A_{k_0})\,dx
  -\lim\limits_{k\to\infty}\sum^N_{i=1}\int\limits_{a_i,\Gamma}f(x)(x-a_i)R(x;A_k)\,dx\bigg \|_p<\epsilon. \eqno (4)$$
Observe that if $1\le i\le N$ and $|\xi|_p=1$, then as in statement (6) of Lemma 2.4, by Lemma 1.3 and (1) imply that 
$x^i_\xi=a_i+\xi\Gamma\not\in D_{\sigma_A}(\epsilon^-)$.
Now let $n_0$ be so large that if $n\ge n_0$ and $p\not\;\mid n$, then
$$\bigg \|\sum\limits^N_{i=1}\bigg\{{1\over n}\sum\limits_{\xi^n=1}f(x^i_\xi)(x^i_\xi-a_i)
R(x^i_\xi;A_{k_0})-
                           \int\limits_{a_i,\Gamma}f(x)(x-a_i)R(x;A_{k_0})\,dx\bigg\}\bigg \|_p<\epsilon.\eqno (5)$$
For  $n\ge n_0$ and  $p\not\;\mid n$, we have
$$\eqalignno{&\sum\limits^N_{i=1}{1\over n}\sum\limits_{\xi^n=1}f(x^i_\xi)(x^i_\xi-a_i)R(x^i_\xi;A)-
\lim\limits_{k\to\infty}\sum^N_{i=1}\int\limits_{a_i,\Gamma}f(x)(x-a_i)R(x;A_k)\,dx &(6)\cr
=&\sum\limits^N_{i=1}{1\over n}\sum\limits_{\xi^n=1}f(x^i_\xi)(x^i_\xi-a_i)
\big\{R(x^i_\xi;A)-R(x^i_\xi;A_{k_0})\big\}+\cr
& \sum\limits^N_{i=1}\bigg\{{1\over n}\sum\limits_{\xi^n=1}f(x^i_\xi)(x^i_\xi-a_i) R(x^i_\xi;A_{k_0})-
                           \int\limits_{a_i,\Gamma}f(x)(x-a_i)R(x;A_{k_0})\,dx\bigg\}+\cr
			   & \sum^N_{i=1}\int\limits_{a_i,\Gamma}f(x)(x-a_i)R(x;A_{k_0})\,dx
			     -\lim\limits_{k\to\infty}\sum^N_{i=1}\int\limits_{a_i,\Gamma}f(x)(x-a_i)R(x;A_k)\,dx.\cr}$$
Thus, for $n\ge n_0$ and  $p\not\;\mid n$, we have (1)-(6) together imply that
$$\bigg \|\sum\limits^N_{i=1}{1\over n}
\sum\limits_{\xi^n=1}f(x^i_\xi)(x^i_\xi-a_i)R(x^i_\xi;A)-
\lim\limits_{k\to\infty}\sum\limits^N_{i=1}\int\limits_{a_i,\Gamma}f(x)(x-a_i)R(x;A_k)\,dx\bigg \|_p<\epsilon.$$
Because $\epsilon>0$ is arbitrary, this completes the proof of the lemma. \vrule height 6pt width 5pt depth 4pt

\medskip
\proclaim 2.7 Definition. Let $\cal A$ be a $p$-adic Banach algebra. Let $A\in{\cal A}$ be $\cal A$-analytic. Let
$(A_k)$ be a sequence in $\cal A$ and $({\cal X}_k)$ a sequence of $p$-adic Banach spaces ${\cal X}_k$ of the form
${\cal X}_k\simeq\Omega_p(J_k)$, with $J_k\not=\emptyset$, such that $\sigma_{A_k}\not=\emptyset$ is compact,
$A_k$ is ${\cal X}_k$-analytic and $\lim\limits_{k\to\infty}\|A-A_k\|_p=0$. Let $A$ have spectrum $\sigma_A$.
Now let $f\in{\cal F}(A)$. Let $0<r\in |\Omega_p|_p$, and let $\Gamma\in\Omega_p$, with $|\Gamma|_p=r$. 
Assume that $a_1,\ldots,a_N$ $\Omega_p$ are given, with
$$\sigma_A\subseteq\bigcup\limits^N_{i=1}D_{a_i}(r), \eqno (1)$$
Assume that $f$ is Krasner analytic on the $D_{a_i}(r)$. By Lemma 2.6, (1) implies that
$$\sum\limits^N_{i=1}\,\int\limits_{a_i,\Gamma}\,f(x)(x-a_i)R(x;A)\,dx=
\lim\limits_{k\to\infty}\sum\limits^N_{i=1}\,\int\limits_{a_i,\Gamma}\,f(x)(x-a_i)R(x;A_k)\,dx.$$
Hence we may define $f(A)$ by
$$f(A)=\sum\limits^N_{i=1}\,\int\limits_{a_i,\Gamma}\,f(x)(x-a_i)R(x;A)\,dx.$$
Lemma 2.5 then implies that this definition of $f(A)$ does not depend on $a_1,\ldots,a_N,r,\Gamma$.

\medskip
\proclaim 2.8 Lemma. Let $\cal A$ be a $p$-adic Banach algebra. Let $A\in {\cal A}$ have spectrum $\sigma_A$,
and assume that $A$ is $\cal A$-analytic. Let $f\in{\cal F}(A)$. Let $0<r$ be in $|\Omega_p|_p$,
and let $\Gamma$ be in $\Omega_p$, 
with $|\Gamma|_p=r$. Assume that $a_1,\ldots,a_N$ in $\Omega_p$ are given, with
$$\sigma_A\subseteq\bigcup\limits^N_{i=1}D_{a_i}(r^-),\eqno (1)$$
where the $D_{a_i}(r)$ are disjoint and $f$ is Krasner analytic on the $D_{a_i}(r)$. Then for any 
$\epsilon >0$ there exists a $\delta >0$ such that if $B\in{\cal A}$ is $\cal A$-analytic with $\|A-B\|_p<\delta$, then
$$\sigma_B\subseteq \bigcup\limits^N_{i=1}D_{a_i}(r),$$
$f\in{\cal F}(A)$ and $\|f(A)-f(B)\|_p<\epsilon$. 

\noindent
{\bf Proof.} Because $f$ is Krasner analytic on each $D_{a_i}(r)$, 
Lemma 1.11 implies that for each $1\le i\le N$, $\max\limits_{x\in D_{a_i}(r)}|f(x)|_p$ is attained when $|x-a_i|_p=r$, 
and hence we may write
$$\rho_i=\max\limits_{x\in D_{a_i}(r)}|f(x)|_p=\max\limits_{|x-a_i|_p=r}|f(x)|_p.$$
Now let $0<\epsilon<r$ be arbitrary. By Lemma 2.1, there exists a positive number $\delta>0$ such that if
$B\in {\cal A}$ and $\|A-B\|_p<\delta$, then
$$\sigma_B\subseteq D_{\sigma_A}(\epsilon^-), \eqno (2)$$
and
$$\|R(x;A)-R(x;B)\|_p<{\epsilon\over r(1+\max\limits_j\rho_j)},\quad x\not\in D_{\sigma_A}(\epsilon^-).\eqno (3)$$
Now, assume that $B\in{\cal A}$ is $\cal A$-analytic, and that $\|A-B\|_p<\delta$. As in (4) of the proof of 
Lemma 2.4, we see that (1) implies 
$$D_{\sigma_A}(r^-)\subseteq \bigcup\limits^N_{i=1}D_{a_i}(r^-). \eqno (4)$$
Then by (2), we have $\sigma_B\subseteq D_{\sigma_A}(\epsilon^-)$, and hence (4) implies that 
$$\sigma_B\subseteq D_{\sigma_A}(\epsilon^-)\subseteq D_{\sigma_A}(r^-)\subseteq
\bigcup\limits^N_{i=1}D_{a_i}(r^-).\eqno (5)$$
Because $B$ is $\cal A$-analytic, $\sigma_B\not=\emptyset$, and hence (5) implies that $f\in{\cal F}(B)$.
Then by (1) and (5), we may use Definition 2.7 to define $f(A)$ and $f(B)$ such that
$$\eqalign{f(A)&=\sum\limits^N_{i=1}\int\limits_{a_i,\Gamma}\,f(x)(x-a_i)R(x;A)\,dx=
                 \lim\limits_{\buildrel{n\to\infty}\over{p\not\;\mid n}}\sum\limits^N_{i=1}{1\over n}
\sum\limits_{\xi^n=1}f(x^i_{\xi})(x^i_{\xi}-a_i)R(x^i_{\xi};A);\cr
             f(B)&=\sum\limits^N_{i=1}\int\limits_{a_i,\Gamma}\,f(x)(x-a_i)R(x;B)\,dx=
	     \lim\limits_{\buildrel{n\to\infty}\over{p\not\;\mid n}}\sum\limits^N_{i=1}{1\over n}
	     \sum\limits_{\xi^n=1}f(x^i_{\xi})(x^i_{\xi}-a_i)R(x^i_{\xi};B).\cr}$$
Then we get 
$$f(A)-f(B)=\lim\limits_{\buildrel{n\to\infty}\over{p\not\;\mid n}}\sum\limits^N_{i=1}{1\over n}
             \sum\limits_{\xi^n=1}f(x^i_{\xi})(x^i_{\xi}-a_i)\big[R(x^i_{\xi};A)-R(x^i_{\xi};B)\big].\eqno(6)$$
Now, we claim that for $1\le i\le N$ and $\xi\in\Omega_p$, with, $|\xi|_p=1$,  
$$x^i_{\xi}\not\in D_{\sigma_A}(\epsilon^-).\eqno (7)$$
The proof of this claim is similar to the proof of statement (6) in Lemma 2.4. 
It then follows from (3), (6) and (7) that 
$$\eqalign{\|f(A)-f(B)\|_p&=\lim\limits_{\buildrel{n\to\infty}\over{p\not\;\mid n}}
                  \bigg\|\sum\limits^N_{i=1}{1\over n}
	          \sum\limits_{\xi^n=1}f(x^i_{\xi})(x^i_{\xi}-a_i)\big[R(x^i_{\xi};A)-R(x^i_{\xi};B)\big]\bigg\|_p\cr
		  &\le\max_{1\le i\le N}\bigg\{(\rho_i r)\bigg[{\epsilon\over r(1+\max\limits_j\rho_j)}\bigg]\bigg\}\cr
		  &<\epsilon.}$$
This completes the proof of the lemma. \vrule height 6pt width 5pt depth 4pt

\medskip 
\proclaim 2.9 Lemma. Suppose that $a_1,\ldots,a_N;b_1,\ldots,b_M$ are in $\Omega_p$, and let $0<r\le s$ be positive
real numbers such that
$$\bigg[\bigcup\limits^N_{i=1}D_{a_i}(r^-)\bigg]\cap\bigg[\bigcup\limits^M_{i=1}D_{b_i}(s^-)\bigg]\not=\emptyset.$$
Then there exists $\emptyset\not=I\subseteq\{\,1,\ldots,N\,\}$ such that 
$$\bigg[\bigcup\limits^N_{i=1}D_{a_i}(r^-)\bigg]\cap\bigg[\bigcup\limits^M_{i=1}D_{b_i}(s^-)\bigg]=
  \bigcup\limits_{i\in I}D_{a_i}(r^-),$$
where for each $i\in I$ there exists a $1\le j\le M$ for which
$$D_{a_i}(r)=D_{a_i}(r)\cap D_{b_j}(s).$$

\noindent
{\bf Proof.} The proof of this lemma is elementary and will be omitted.  \vrule height 6pt width 5pt depth 4pt

\medskip
\proclaim 2.10 Extended Spectral Theorem. Let $cal A$ be a $p$-adic Banach algebra. Let $A\in {\cal A}$ be 
$\cal A$-analytic, with spectrum $\sigma_A$. Let $f\in {\cal F}(\sigma_A)$. 
Let $0<r$ be in $|\Omega_p|_p$,
and let $\Gamma$ be in $\Omega_p$, with $|\Gamma|_p=r$. Assume that $a_1,\ldots,a_N$ in $\Omega_p$ are given, with
$$\sigma_A\subseteq\bigcup\limits^N_{i=1}D_{a_i}(r^-),\eqno (1)$$
where the $D_{a_i}(r)$ are disjoint and $f$ is Krasner analytic on each $D_{a_i}(r)$. By Definition 2.7, define
$$f(A)=\sum\limits^N_{i=1}\int\limits_{a_i,\Gamma}\,f(x)(x-a_i)R(x;A)\,dx.$$
Let $g\in {\cal F}(A)$, and let $\alpha,\beta \in \Omega_p$. Then $\alpha f+\beta g, f\cdot g\in {\cal F}(A)$, and
\item\item{(2)} $(\alpha f+\beta g)(A)=\alpha f(A)+\beta g(A)$;
\item\item{(3)} $(f\cdot g)(A)=f(A)g(A)$;
\item\item{(4)} The mapping $h\mapsto h(A)$ is continuous on $B_r(\sigma_A)$.

\noindent
{\bf Proof.} Let $g\in{\cal F}(A)$. Let $0<s$ be in $|\Omega_p|_p$, and suppose that $b_1,\ldots,b_M$ in $\Omega_p$
are such that the $D_{b_i}(s)$ are disjoint, $g$ is Krasner analytic on each $D_{b_i}(s)$, and 
$$\sigma_A\subseteq \bigcup\limits^M_{i=1}D_{b_i}(s^-). \eqno (5)$$
Without loss of generality we may assume that $r\le s$. Because $A$ is $\cal A$-analytic, $\sigma_A\not=\emptyset$. 
Hence, by Lemma 2.9, (1) 
and (5) together imply that there exists $\emptyset\not=J\subseteq\{\,1,\ldots,N\,\}$ such that
$$\bigg[\bigcup\limits^N_{i=1}D_{a_i}(r^-)\bigg]\cap\bigg[\bigcup\limits^M_{i=1}D_{b_i}(s^-)\bigg]=
  \bigcup\limits_{i\in J}D_{a_i}(r^-),$$
where for each $i\in I$ there exists a $1\le j\le M$ for which
$$D_{a_i}(r)=D_{a_i}(r)\cap D_{b_j}(s).$$
Therefore, if $i\in J$, then $g$ is Krasner analytic on $D_{a_i}(r)$. Hence we may assume
without loss of generality that $g$ is Krasner analytic on each $D_{a_i}(r)$. Finally, by Lemma 1.3 and (1), 
we may assume without loss of generality that   
$$\sigma_A\subseteq D_{\sigma_A}(r^-)=\bigcup\limits^N_{i=1}D_{a_i}(r^-)\hbox{ and }
D_{\sigma_A}(r)=\bigcup\limits^N_{i=1}D_{a_i}(r). \eqno (6)$$
Now let $\alpha,\beta\in \Omega_p$, then $\alpha f+\beta g$ and $f\cdot g$ are Krasner analytic on the 
$D_{a_i}(r)$. Because $A$ is $\cal A$-analytic, there exists sequence $(A_k)$ in $\cal A$ and $({\cal X}_k)$ a sequence of
$p$-adic Banach spaces ${\cal X}_k$ of the form ${\cal X}_k\simeq\Omega_p(J_k)$, with $J_k\not=\emptyset$, such that
such that $\sigma_{A_k}\not=\emptyset$ is compact, $A_k$ is ${\cal X}_k$-analytic and 
$\lim\limits_{k\to\infty}\|B-B_k\|_p=0$. 
Then by Theorem 2.1 and (6), we may assume that for all $k$, 
$$\sigma_{A_k}\subseteq D_{\sigma_A}(r^-)\subseteq \bigcup\limits^N_{i=1}D_{a_i}(r^-). \eqno (7)$$
Then by Lemma 1.3, (7) implies that for all $k$ there exist $\emptyset\not= I_k$ such that
$I_k\subseteq \{\,1,\ldots,N\,\}$ and
$$\sigma_{A_k}\subseteq D_{\sigma_{A_k}} (r^-)=\bigcup\limits_{i\in I_k}D_{a_i}(r^-)\hbox{ and }
D_{\sigma_{A_k}} (r)=\bigcup\limits_{i\in I_k}D_{a_i}(r). \eqno  (8)$$
Then by Lemma 2.4 and Theorem 1.21 (Vishik's Spectral Theorem II), (8) implies that if $h$ is any $\Omega_p$-valued 
function that is Krasner analytic on each of the $D_{a_i}(r)$, then
$$\eqalignno{h(A)&=\lim\limits_{k\to\infty}\sum\limits^N_{i=1}\,\int\limits_{a_i,\Gamma}\,h(x)(x-a_i)R(x;A_k)\,dx
              &(9)\cr
   &=\lim\limits_{k\to\infty}\sum\limits^{N_k}_{i=1}\,\int\limits_{a_{k,i},\Gamma}\,h(x)(x-a_{k,i})R(x;A_k)\,dx.\cr
   &=\lim\limits_{k\to\infty}h(A_k).\cr}$$
Then (9) gives
$$\eqalign{(\alpha f+\beta g)(A)&=\lim\limits_{k\to\infty}(\alpha f+\beta g)(A_k)\cr
                                &=\lim\limits_{k\to\infty}(\alpha f(A_k)+\beta g(A_k))\cr
				&=\alpha \lim\limits_{k\to\infty} f(A_k)+\beta \lim\limits_{k\to\infty} g(A_k)\cr
				&=\alpha f(A)+\beta g(A).\cr}$$
This proves (2). By Theorem 1.21, (5) implies that
$$\eqalign{(f\cdot g)(A)&=\lim\limits_{k\to\infty}(f\cdot g)(A_k)\cr
                                &=\lim\limits_{k\to\infty}f(A_k) g(A_k)\cr
		                &=f(A) g(A).\cr}$$
This proves (3). To prove (4), let $h\in B_r(\sigma_p)$. 
Then by (6), $h$ is Krasner analytic on each $D_{a_i}(r)$, $i\in I$. Lemma 1.11 then implies 
that for each $j\in I$, $\max\limits_{x\in D_{a_j}(r)}|h(x)|_p$ is attained when 
$|x-a_j|_p=r$, and hence we may write
$$\rho_j(h)=\max\limits_{x\in D_{a_j}(r)}|h(x)|_p=\max\limits_{|x-a_j|_p=r}|h(x)|_p \hbox{ and }
\|h\|_r=\max\limits_{i\in I}\rho_i(h).$$
By Theorem 2.1, there exists a constant $N_r>0$ such that $\|R(x;A)\|_p\le N_r$ for all $x\not\in D_{\sigma_A}(r^-)$.
Now let $\epsilon > 0$ be given, and define $\displaystyle{\delta={\epsilon\over r N_r}}$. 
Let $h_0,h\in B_r(\sigma_A)$,
with $\|h-h_0\|_r<\delta$. Let $\Gamma\in\Omega_p$, with $|\Gamma|_p=r$. For $1\le i\le N$ and $\xi\in\Omega_p$, with
$|\xi|_p=1$, set $x^i_\xi=a_i+\xi\Gamma$; then by (6),  $x^i_\xi\not\in D_{\sigma_A}(r^-)$, hence
by Lemma 2.6, we have
$$\eqalign{&\|h(A)-h_0(A)\|_p\cr
=&\|(h-h_0)(A)\|_p\cr
=&     \lim\limits_{\buildrel{n\to\infty}\over{p\not\;\mid n}}\bigg\|\sum\limits^N_{i=1}{1\over n}
     \sum\limits_{\xi^n=1}[h(x^i_{\xi})-h_0(x^i_{\xi})](x^i_{\xi}-a_i)R(x^i_{\xi};A)\bigg\|_p\cr
=&\lim\limits_{\buildrel{n\to\infty}\over{p\not\;\mid n}}\bigg\|\sum\limits^N_{i=1}{1\over n}
     \sum\limits_{\xi^n=1}[h(x^i_{\xi})-h_0(x^i_{\xi})](\xi\Gamma)R(x^i_{\xi};A)\bigg\|_p\cr
\le &\|h-h_0\|_r(r N_r)\cr
<&{\epsilon\over r N_r}(r N_r)=\epsilon.\cr}$$
This proves (4). \vrule height 6pt width 5pt depth 4pt
\vskip 1cm
\centerline{\bf 3. $p$-ADIC UHF AND TUHF BANACH ALGEBRAS}
\vskip 1cm
\proclaim 3.1 Definition. A {\bf $p$-adic UHF algebra} over $\Omega_p$ is a unital $p$-adic Banach algebra of the form
$${\cal A}=\overline{\bigcup\limits^\infty_{n=1}{\cal M}_n},$$
where each $({\cal M}_n)$ is an increasing sequence of $p$-adic Banach subalgebras of $\cal A$, such that each
${\cal M}_n$ contains the identity of $\cal A$ and is algebraically isomorphic as 
an $\Omega_p$-algebra to $M_{p_n}(\Omega_p)$. A {\bf $p$-adic TUHF algebra} over 
$\Omega_p$ is a unital $p$-adic Banach algebra of the form
$${\cal T}=\overline{\bigcup\limits^\infty_{n=1}{\cal T}_n},$$
where each $({\cal T}_n)$ is an increasing sequence of $p$-adic Banach subalgebras of $\cal T$, such that each
${\cal T}_n$ contains the identity of $\cal T$ and is algebraically isomorphic as an 
$\Omega_p$-algebra to $T_{p_n}(\Omega_p)$.

\medskip
\proclaim 3.2 Lemma. Let $\cal X$ be a $p$-adic Banach space. with  ${\cal X}\simeq\Omega_p(J)$, where $J\not=\emptyset$
is a finite set. $\cal A$ be a $p$-adic Banach algebra and assume that there exists an $\Omega_p$-monomorphism
$\theta:{\cal A}\to B({\cal X})$. Then any $A\in {\cal A}$ is ${\cal X}$-analytic.

\medskip
\noindent
{\bf Proof.} Let $A\in{\cal A}$. Because $B({\cal X})$ is finite-dimensional, the results of finite-dimensional
spectral theory apply to operators in $B({\cal X})$. Therefore there exist $x_1,\ldots,x_N$ in $\Omega_p$ such that
$\sigma_A=\{\,x_1,\ldots,x_N\,\}$, where $\sigma_A$ is the spectrum of $A$. Moreover, there exist positive integers
$\nu(x_1),\ldots,\nu(x_N)$ and idempotents $E(x_1),\ldots,E(x_N)$ in  $B({\cal X})$ such that for $x\not\in\sigma_A$,
$$R(x;\theta(A))=\sum\limits^N_{j=1}\sum\limits^{\nu(x_j)-1}_{\nu=0}
{(\theta(A)-x_j)^\nu\over(x-x_j)^{\nu+1}}E(x_j).\eqno (1)$$
For a proof of (1), see {\bf [DS]}: Theorem 10, p. 560. Now fix $y\in{\cal X}$ and $h\in{\cal X}^*$. Then for
$x\not\in\sigma_A$, (1) implies that
$$h(R(x;\theta(A))y)=
\sum\limits^N_{j=1}\sum\limits^{\nu(x_j)-1}_{\nu=0} {h_{j\nu}\over(x-x_j)^{\nu+1}},\eqno (2)$$
where for $1\le j\le N$ and $0\le \nu\le \nu(x_j)-1$, $h_{j\nu}=h((\theta(A)-x_j)^\nu E(x_j)y)$. Now, each of the
functions $\displaystyle{{h_{j\nu}\over(x-x_j)^{\nu+1}}}$ in (2) is clearly in $H_0(\overline{\sigma}_{\theta(A)})$,
hence (2) implies that the function $R_{h,y}(x)=h(R(x;\theta(A))y)$ belongs to $H_0(\overline{\sigma}_{\theta(A)})$.
Because $y\in{\cal X}$ and $h\in{\cal X}^*$ are arbitrary, we see that $\theta(A)$ is analytic, and hence $A$ is
${\cal X}$-analytic. \vrule height 6pt width 5pt depth 4pt

\medskip
\proclaim 3.3 Spectral Theorem for $p$-adic UHF and TUHF Algebras. Let $\cal A$ be a $p$-adic {\rm UHF}
or a $p$-adic {\rm TUHF} Banach algebra.
Let $A\in{\cal A}$. Then $A$ is $\cal A$-analytic with spectrum $\sigma_A\not=\emptyset$. Let $f\in {\cal F}(\sigma_A)$.
Let $0<r$ be in $|\Omega_p|_p$, and let $\Gamma$ be in $\Omega_p$, with $|\Gamma|_p=r$. Assume that 
$a_1,\ldots,a_N$ in $\Omega_p$ are given, with
$$\sigma_A\subseteq\bigcup\limits^N_{i=1}D_{a_i}(r^-),\eqno (1)$$
where the $D_{a_i}(r)$ are disjoint and $f$ is Krasner analytic on each $D_{a_i}(r)$. By Definition 2.7, define
$$f(A)=\sum\limits^N_{i=1}\int\limits_{a_i,\Gamma}\,f(x)(x-a_i)R(x;A)\,dx.$$
Let $g\in {\cal F}(A)$, and let $\alpha,\beta \in \Omega_p$. Then $\alpha f+\beta g, f\cdot g\in {\cal F}(A)$, and
\item\item{(2)} $(\alpha f+\beta g)(A)=\alpha f(A)+\beta g(A)$;
\item\item{(3)} $(f\cdot g)(A)=f(A)g(A)$;
\item\item{(4)} The mapping $h\mapsto h(A)$ is continuous on $B_r(\sigma_A)$.

\medskip
\noindent
{\bf Proof.} $\cal A$ has the from
$${\cal A}=\overline{\bigcup\limits^\infty_{n=1}{\cal A}_n},$$
where for each $n$, there exists a finite nonempty set $J_n$ and an $\Omega_p$-monomorphism
$\theta_n:{\cal A}_n\to B({\cal X}_n)$, with ${\cal X}_n\simeq\Omega_p(J_n)$; because $B({\cal X}_n)$ is 
finite-dimensional, $\theta_n$ is bicontinuous. It follows form Lemma 3.2 that for all
$n$, the members of ${\cal A}_n$ are ${\cal X}_n$-analytic. Let $A\in{\cal A}$, there exists a sequence $(A_n)$ in
$$\bigcup\limits^\infty_{n=1}{\cal A}_n$$
such that $\lim\limits_{n\to\infty}\|A-A_n\|_p=0$, therefore, by Definition 1.19, $A$ is $\cal A$-analytic. Properties
(2)-(4) are then consequences of applying Theorem 2.10 to $A$. \vrule height 6pt width 5pt depth 4pt

\medskip
\proclaim 3.4 Theorem. Let $\cal A$ be a $p$-adic {\rm UHF} or a $p$-adic {\rm TUHF} Banach algebra. Let $A\in {\cal A}$, 
then $A$ is $\cal A$-analytic
with with spectrum $\sigma_A$. Let $f\in{\cal F}(A)$. Let $0<r$ be in $|\Omega_p|_p$, and let $\Gamma$ be in $\Omega_p$,
with $|\Gamma|_p=r$. Assume that $a_1,\ldots,a_N$ in $\Omega_p$ are given, with
$$\sigma_A\subseteq\bigcup\limits^N_{i=1}D_{a_i}(r^-),\eqno (1)$$
where the $D_{a_i}(r)$ are disjoint and $f$ is Krasner analytic on the $D_{a_i}(r)$. Then for any
$\epsilon >0$ there exists a $\delta >0$ such that if $B\in{\cal A}$, with $\|A-B\|_p<\delta$, then
$$\sigma_B\subseteq \bigcup\limits^N_{i=1}D_{a_i}(r),$$
$f\in{\cal F}(B)$ and $\|f(A)-f(B)\|_p<\epsilon$.

\medskip
\noindent
{\bf Proof.} In the proof of Theorem 3.3 it was demonstrated that every member of $\cal A$ is $\cal A$-analytic.
Therefore the present lemma follows from an application of Lemma 2.8.  \vrule height 6pt width 5pt depth 4pt

As mentioned in the introduction of the present article, Theorem 3.3 and Theorem 3.4 can be used to transfer, 
{\it mutatis mutandis}, the results of {\bf [B1]} to $p$-adic TUHF algebras, which is the content of {\bf [B2]}.

\bigskip
\medskip
\centerline{REFERENCES}
\medskip
\parindent=.5cm
\item{\bf [B1].} R. L. Baker, {\it On Certain Banach Limits of Triangular Matrix Algebras}, Houston J. Math., {\bf 23}, 
No.1, (1997), 127-141.
\item{\bf [B2].} R. L. Baker, {\it The Classification of $p$-adic TUHF Banach Algebras}, preprint in preparation.
\item{\bf [DS].} N. Dunford and J. T. Schwartz, {\it Linear Operators Part I}, Interscience Publishers, Inc., New York.
\item{\bf [K].} N. Koblitz, {\it $p$-adic Analysis: a Short Course on Recent Work}, London Mathematical Society
Lecture Note Series 46, Cambridge University Press. 
\item{\bf [NBB].} L. Narici, E. Beckenstein, and G. Bachman, {\it Functional Analysis and Valuation Theory}, Marcel 
Dekker, Inc., New York 1971.
\bigskip
\parindent .5cm
Department of Mathematics, University of Iowa, Iowa City, Iowa 52242

{\it E-mail address}: {\bf baker@math.uiowa.edu}

\vfill\eject

\end